\documentclass[12pt]{amsart}
\usepackage{amssymb,amscd,amsxtra,calc}
\usepackage{cmmib57}
\usepackage{graphicx}
\usepackage{amssymb,amsmath,amsfonts,mathrsfs,enumerate,xcolor}
\usepackage{amsthm}

\usepackage{soul}
\usepackage[citecolor=Navy]{hyperref}
\usepackage[ruled,vlined]{algorithm2e}

\usepackage{lipsum}
\usepackage{amsfonts}
\usepackage{graphicx}
\usepackage{epstopdf}
\usepackage{algorithmic}
\usepackage{float}
\restylefloat{table}
\usepackage{placeins}
\usepackage{array}
\ifpdf
  \DeclareGraphicsExtensions{.eps,.pdf,.png,.jpg}
\else
  \DeclareGraphicsExtensions{.eps}
\fi

\usepackage{subcaption}

\theoremstyle{plain}
    \newtheorem{thm}{Theorem}[section]

    \newtheorem{lemma}[thm]{Lemma}
    \newtheorem{proposition}[thm]{Proposition}
    \newtheorem{conjecture}[thm]{Conjecture}

\theoremstyle{definition}

    \newtheorem{remark}[thm]{Remark}
\theoremstyle{remark}

\newcommand{\authorfootnotes}{\renewcommand\thefootnote{\@fnsymbol\c@footnote}}

\usepackage{caption} 
\title[]{Some interesting birational morphisms of smooth affine quadric $3$-folds}

\author{Cinzia Bisi} 
\address{Department of Mathematics and Computer Science, University of Ferrara, Italy}
\email{bsicnz@unife.it}
\author{Jonathan D. Hauenstein} 
\address{Department of Applied and Computational Mathematics and Statistics, University of Notre Dame, USA}
\email{hauenstein@nd.edu}
\author{Tuyen Trung Truong}
\address{Department of Mathematics, University of Oslo, Norway}
\email{tuyentt@math.uio.no}
\thanks{}
\date{\today}
\begin{document}
\maketitle
{\centering\footnotesize To the late Nessim Sibony\par}
\begin{abstract} We study a family of birational maps of smooth affine quadric 3-folds, {over the complex numbers}, of the form $x_1x_4-x_2x_3=$ constant, which seems to have some (among many others) interesting/unexpected characters:  a) they are cohomologically hyperbolic, b) their second dynamical degree is an algebraic number but not an algebraic integer, and c) the logarithmic growth of their periodic points is strictly smaller than their algebraic entropy. These maps are restrictions of a polynomial map on $\mathbb{C}^4$ preserving each of the quadrics. The study in this paper is a mixture of rigorous and experimental ones, where for the experimental study we rely on Bertini which is a reliable and fast software for expensive numerical calculations in complex algebraic geometry.  
\end{abstract}











\section{Introduction}
A main theme in Complex Dynamics is that of equidistribution for periodic points. Roughly speaking, it is expected that if a rational map (more generally, meromorphic map) $f:X\dashrightarrow X$ of a complex projective manifold (more generally, compact K\"ahler manifold) satisfies a certain type of hyperbolicity, then it has an equilibrium measure with no mass on proper Zariski subsets, and to which the average of Dirac measures on hyperbolic periodic points converges. In general, one expects that the logarithmic growth of (hyperbolic) isolated periodic points of a map is a lower bound for its topological entropy  (a very fundamental dynamical invariant). These are characteristics of the complexity of a map. 

To make precise the above statement, we first recall the definition of dynamical degrees of a map. Let $X$ be a compact K\"ahler manifold of dimension $m$, and $f:X\dashrightarrow X$ a dominant meromorphic map (i.e.: there is a proper Zariski subset $I(f)\subset X$ - called the indeterminacy set of $f$ - so that $f$ is given as a holomorphic map $f:X\backslash I(f)\rightarrow X$ and the image of $X\backslash I(f)$ is dense in $X$). Let $\omega$ be the cohomological class of a K\"ahler form on $X$. For a number $0\leq j\leq m$, the following limit exists and is independent of the choice of the K\"ahler form (see \cite{russakovskii-shiffman} when $X=$ a projective space, \cite{dinh-sibony10} and \cite{dinh-sibony1} for the general case)
\begin{eqnarray*}
\lim _{n\rightarrow \infty}||(f^n)^*(\omega ^j)||_{H^{2j}(X,\mathbb{C})}^{1/n}. 
\end{eqnarray*}
Here, $f^n=f\circ f\circ \ldots \circ f$ (n times) is the n-th iterate of $f$, and $||.||$ is any norm on the finite dimensional vector space $H^{2j}(X,\mathbb{C})$. We call the above limit the $j$-th dynamical degree of $f$, and denote by $\lambda _j(f)$. In the case $X$ is a complex projective manifold, one can choose $\omega$ to be the class of an ample divisor and replace $H^{2j}(X,\mathbb{C})$ by $N_{\mathbb{R}}^j(X)$ the space of algebraic cycles of codimension $j$ (with real coefficients)  modulo the numerical equivalence. 
{One always has} $\lambda_0(f)=1$ while $\lambda _m(f)$ is the topological degree of $f$ (and hence is $1$ when $f$ is bimeromorphic). It is known that $\lambda _j(f)$'s are bimeromorphic invariants (\cite{dinh-sibony10} and \cite{dinh-sibony1}), that is if $\pi :Y\dashrightarrow X$ is a bimeromorphic map, and $f_Y$ is the lifting of $f$ to $Y$ (i.e. $f_Y=\pi ^{-1}\circ f\circ \pi$), then $\lambda _j(f_Y)=\lambda _j(f)$ for all $j$.  There are also arithmetic analogs of these dynamical degrees, where the above mentioned results are largely unknown, and there are many challenging conjectures around (see \cite{kawaguchi-silverman}\cite{dang-etal}). Dynamical degrees can also be defined for maps on fields of positive characteristics \cite{truong}\cite{dang}, and the version for correspondence in \cite{truong} thus can be used to provide a generalization of Weil's Riemann hypothesis \cite{truong1}\cite{hu-truong1}\cite{hu-truong2}. Dynamical degrees provide a useful way to compute topological entropy of holomorphic selfmaps of compact K\"ahler manifolds \cite{gromov}\cite{yomdin}. Inspired by this fact, for a general meromorphic map, the logarithm of the maximum of its dynamical degrees is named its {\bf algebraic entropy} in \cite{bellon-etal2}, and this is widely practiced in the literature. Moreover, they provide an upper bound for the topological entropy of dominant meromorphic maps of compact K\"ahler manifolds \cite{dinh-sibony10}\cite{dinh-sibony1}, and the same upper bound holds in the Berkovich space setting \cite{favre-rivera}\cite{favre-etal}. Besides being of interest in pure mathematics, birational maps appear naturally in some physical models (in lattice statistical mechanics), and their dynamical degrees are an indication of the complexity of these models, see e.g. \cite{bellon-etal}\cite{bellon-etal2}\cite{bellon-viallet}\cite{bedford-kim}\cite{bedford-truong}\cite{truong3}.    

A dominant meromorphic map $f:X\dashrightarrow X$ is cohomologically hyperbolic if there is one index $j$ so that $\lambda _j(f)>\max _{i\not= j}\lambda _i(f)$. Cohomological hyperbolicity is a cohomological version of the well known notion of hyperbolic dynamics in differentiable dynamical systems, where many results around periodic points, equilibrium measures and topological entropy are known, see a survey in \cite{hasselblatt-pesin}. The rough idea is that since algebraic dynamical systems are more rigid than smooth dynamical systems,  hyperbolicity of an algebraic dynamical system may be detected by the easier invariants on cohomology groups.   

With this preparation, an explicit statement of a major folklore conjecture, which attracts a lot of attention and work, is the following:

\begin{conjecture}[Folklore Conjecture] Let $X$ be a compact K\"ahler manifold of dimension $m$, and $f:X\dashrightarrow X$ a dominant meromorphic map.  Assume that $f$ is cohomologically hyperbolic. Then the following are true: 

1) There is a probability measure $\mu$ on $X$ with no mass on proper analytic subsets (hence, one can push it forward by $f$) which is invariant by $f$, i.e. $f_*\mu =\mu$. 

2) Let $HP_n(f)$ be the set of hyperbolic periodic points of period $n$ of $f$. Then: 

a) The exponential growth of $\sharp HP_n(f)$  is $\lambda (f)$, where $\lambda (f)=\max _j\lambda _j(f)$. This means that
\begin{eqnarray*}
\lim _{n\rightarrow\infty}\frac{\log \sharp HP_n(f)}{n}=\log \lambda (f).
\end{eqnarray*}

b) The hyperbolic periodic points of $f$ equidistribute to $\mu$. That is, 
\begin{eqnarray*}
\lim _{n\rightarrow\infty}\frac{1}{\sharp HP_n(f)}\sum _{x\in HP_n(f)}\delta _{x}=\mu.
\end{eqnarray*}
Here, $\delta _x$ is the Dirac measure at $x$. 

c) The number of other isolated periodic points of $f$ of order $n$ is negligible compared to that of $\sharp HP_n(f)$. 

\label{Conjecture1}\end{conjecture}

A stronger version (which also involves the topological entropy and Lyapunov exponents) of this conjecture was stated in \cite{guedj}. A few representatives from the large known results in the literature resolving Conjecture \ref{Conjecture1} in the affirmative are: dimension 1 (see \cite{brolin}\cite{lyubich}\cite{freie-etal}); H\'enon maps in dimension 2 (see \cite{bedford-etal1}\cite{bedford-etal2}); automorphisms of K3 surfaces (see \cite{cantat}); birational maps of surfaces satisfying a certain condition on Green currents (see \cite{dujardin}, which relies on a precise estimate for the number of isolated periodic points: the lower bound is obtained therein by using laminar currents, and the upper bound is later provided by \cite{iwasaki-uehara}\cite{dinh-etal2}); a large class of maps of surfaces satisfying an energy condition (see \cite{diller-etal}); meromorphic maps of compact K\"ahler manifolds for which $\lambda _m>\lambda _{m-1}$ (see \cite{dinh-etal1}); and analogs of H\'enon-maps in higher dimensions (see \cite{dinh-sibony3}). For a comprehensive survey, see \cite{dinh-sibony2}, where other topics of equidistribution - besides that of periodic points - are also discussed. Despite all of these partial results and many efforts, Conjecture \ref{Conjecture1} is still largely open. 

There are some indications in the literature that Conjecture \ref{Conjecture1} may be too strong to be true. For example, one related conjecture is that $\lambda _1(f)$ is an algebraic number \cite{bellon-viallet}. However, recent work \cite{bell-etal}\cite{bell-etal2} shows the existence of maps (can be chosen to be birational maps of $\mathbb{P}^d$ where $d\geq 3$) for which $\lambda _1(f)$ is a transcendental number, and hence the mentioned related conjecture does not hold. On the other hand, it is not known if the mentioned maps in \cite{bell-etal}\cite{bell-etal2} provide counter-examples to Conjecture \ref{Conjecture1}, because these papers treat only the first dynamical degrees and concern neither periodic points nor equilibrium measures. 

It is noteworthy that for the various maps in the literature where dynamical degrees can be actually computed (for a birational maps of surfaces a general procedure using point blowups has been given in \cite{diller-favre}, but besides that one must in general treat case by case), the dynamical degrees are either {\bf algebraic integers} (i.e. roots of a polynomial $p(t)$ whose coefficients are integers and the leading coefficient is $1$) or  {\bf transcendental numbers} (i.e. not a root of any polynomial with integer coefficients). For some special maps, there are even particular speculations about that the dynamical degrees are algebraic integers of special types. Therefore, it is natural to ask what is the situation for the numbers in between these two types, that is algebraic numbers which are not algebraic integers. For example, recall that a real number $\lambda$ is a weak Perron number if it is an algebraic integer and any of its Galois conjugate $\mu$ satisfies $|\mu |\leq \lambda$. There is the following conjecture \cite{blanc-Santen2} (for the first dynamical degree) and  \cite{dang-favre} (for all dynamical degrees), concerning the very actively studied polynomial maps of affine spaces: 
\begin{conjecture}
If $f:\mathbb{C}^d\rightarrow \mathbb{C}^d$ is a polynomial map, then the dynamical degrees (of the extension of $f$ to $\mathbb{P}^d$) are all weak Perron numbers (in particular, algebraic integers). 
\label{ConjectureAlgebraicInteger}\end{conjecture}

For the first dynamical degree, the above conjecture holds in dimension 2 (see \cite{friedland-milnor}\cite{favre-jonsson}) and some automorphisms in $\mathbb{C}^3$ (see \cite{maegawa}), and some special automorphisms in higher dimensions (\cite{blanc-Santen}\cite{blanc-Santen2} and an unpublished result by Mattias Jonsson cited in these papers). For proper polynomials of $\mathbb{C}^m$ so that $\lambda _1^2>\lambda _2$, \cite{dang-favre} shows that $\lambda _1(f)$ is an algebraic number of degree at most $m$.  For all dynamical degrees, one  non-trivial  case  is that of a monomial map $f$ (see \cite{favre-wulcan} \cite{lin})  where all of its dynamical degrees are known to be the product of absolute values of some eigenvalues of the integer matrix defining $f$, and hence are all algebraic integers (recall that algebraic integers form a ring). 

The purpose of this paper is to present a  family of birational maps of smooth quadric 3-folds which are a candidate counter-example for Conjecture \ref{Conjecture1}. At the same time, these birational maps come from polynomial maps on $\mathbb{C}^4$ whose second and third dynamical degrees seem to be an algebraic number but not an algebraic integer (and hence, can also be a counter-example for Conjecture \ref{ConjectureAlgebraicInteger}). As such, the mentioned methods in the last part of the previous paragraph cannot be used to compute the conjectured second dynamical degree, and new ideas have to be developed for this task. 

A disadvantage of the statement of Conjecture \ref{Conjecture1} is that some parts of it depend on an unspecified birational map $\pi :Y\dashrightarrow X$ and an unspecified equilibrium measure $\mu$. As will be seen in Section 3, the following conjecture is a consequence of Conjecture \ref{Conjecture1}, and concerns only the exponential growth of the set of isolated periodic points (which in theory can be explicitly found) and is independent of the birational model of a given map. 

\begin{conjecture} Let $X$ be a smooth complex projective variety  of dimension $d$, and let $f:X\dashrightarrow X$ be a dominant rational map. Assume that $f$ is cohomologically hyperbolic. Define $\lambda (f)=\max _{j}\lambda _j(f)$ and $IsoPer_n(f)$ the set of isolated periodic points of period $n$ of $f$ (multiplicities accounted). Let $Z\subset X$ be any Zariski open dense set. Then 
\begin{eqnarray*}
\limsup _{n\rightarrow\infty}\frac{\log \sharp (IsoPer_n(f)\cap Z)}{n}=\log \lambda (f).
\end{eqnarray*}
\label{Conjecture2}\end{conjecture}

Note that one part of the above conjecture is known \cite{dinh-etal2} (valid in the K\"ahler setting, and the proof therein uses the theory of tangent currents developed in \cite{dinh-sibony4}): we always have
\begin{eqnarray*}
\limsup _{n\rightarrow\infty}\frac{\log \sharp (IsoPer_n(f)\cap Z)}{n}\leq \log \lambda (f).
\end{eqnarray*}
(On the other hand, on differentiable dynamical systems, there are no such upper bounds on the number of periodic points, see \cite{kaloshin}.)

Hence, to disprove Conjecture \ref{Conjecture1} it suffices to display a counter-example to Conjecture \ref{Conjecture2}. We next present our candidate counter-examples for Conjectures \ref{Conjecture2} (and hence \ref{Conjecture1}) and \ref{ConjectureAlgebraicInteger}, which are strikingly simple. We start with a polynomial $F:\mathbb{C}^4\rightarrow \mathbb{C}^4$  given by: 
$$F(x_1,x_2,x_3,x_4)=\left(x_2,-x_4, x_1-x_1x_2^2,-x_3+x_1x_2x_4\right).$$ This is a birational map, with inverse $$F^{-1}(x_1,x_2,x_3,x_4)=\left(\frac{-x_3}{-1+x_1^2},x_1,\frac{x_1x_2x_3+x_4-x_1^2x_4}{-1+x_1^2},-x_2\right).$$

\begin{remark} The above map $F$ is an element of a more general family: 
\begin{equation}
G(x_1,x_2,x_3,x_4)=(x_2,-x_4,x_1-x_2P(x_1,x_2,x_3,x_4),-x_3+x_4P(x_1,x_2,x_3,x_4)).
\label{EquationBirational}\end{equation}

Here $P(x_1,x_2,x_3,x_4)$ is a polynomial (or more generally a rational function) of the form 
$$P(x_1,x_2,x_3,x_4)=x_1Q_1(x_2,x_4)+x_3Q_3(x_2,x_4)+R(x_2,x_4).$$

Some special automorphisms of $\mathbb{C}^4$ in the family $G$ have been studied in \cite{bisi-etal}. 
\label{Remark1}\end{remark}

Let $\phi : \mathbb{C}^4\rightarrow \mathbb{C}$ be the map $\phi (x_1,x_2,x_3,x_4)=x_1x_4-x_2x_3$. It can be checked that the fibres of $\phi$ are invariant by $f$, that is $\phi =\phi \circ F$. In other words, via $\phi$, $F$ is semi-conjugate to the identity map on $\mathbb{C}$.   For a generic $c\in \mathbb{C}$, let $Z_c=\phi ^{-1}(c)=\{x_1x_4-x_2x_3=c\}$ $\subset \mathbb{C}^4$ be the fibre of $\phi$ over $c$. Let $X_c=$ the closure in $\mathbb{P}^4$ of $Z_c$ (hence, in homogeneous coordinates $[x_1:x_2:x_3:x_4:z]$ of $\mathbb{P}^4$, $X_c$ is given by $x_1x_4-x_2x_3=cz^2$). Then both $X_c$ and $Z_c$ are smooth, invariant by $F$, and $Z_c$ is a Zariski open dense subset of $X_c$. The remaining of this paper is to study the following conjecture for the map $f_c=F|_{Z_c}:Z_c\rightarrow Z_c$. Note that since $f_c$ is a birational map in dimension 3, we have $\lambda _0(f_c)=\lambda _3(f_c)=1$. Hence, the only non-trivial dynamical degrees of $f_c$ are $\lambda _1(f_c), \lambda _2(f_c)$. 

\begin{conjecture} Let $f_c$ be the above map, for a generic value $c\in \mathbb{C}$, and $\widehat{f}_c$ its extension to $X_c$. Then we have: 

1) The first dynamical degree $\lambda _1(\widehat{f}_c)$ is the largest root $\zeta _1$ of the polynomial $t^3-t^2-t-1$, and $\zeta _1$ is approximately $1.8393$. 

2) The second dynamical degree $\lambda _2(\widehat{f}_c)$ is the largest root $\zeta _2$ of the polynomial $2t^3-3(t^2-1)-4t$. This polynomial is irreducible over $\mathbb{Q}$, and hence $\lambda _2(f_c)$ is an algebraic number but not an algebraic integer. Here, $\zeta _2$ is approximately $2.1108$. 

3) Moreover, 

{\bf Stronger estimate:}

$$\limsup _{n\rightarrow\infty}\frac{ \log (\sharp IsoFix_n(f_c))}{n}\leq \log 2.108.$$

{\bf Weaker estimate:}

$$\limsup _{n\rightarrow\infty}\frac{ \log (\sharp IsoFix_{2n+1}(f_c))}{2n+1}\leq \log 2.108.$$

\label{Conjecture3}\end{conjecture}

Part 1 of the above conjecture is solved in Lemma \ref{Lemma4}. If $f_c$ indeed satisfies Conjecture \ref{Conjecture3}, then it is a primitive map, see Section 3 for detail. Conjecture \ref{Conjecture3}, if holds,  implies that Conjecture \ref{Conjecture2} (and hence also Conjecture \ref{Conjecture1})  and Conjecture \ref{ConjectureAlgebraicInteger} do not hold, see Section 3. Indeed, the experimental results in Section 3 seem to indicate that the limit $$\lim _{n\rightarrow\infty}\frac{ \log (\sharp IsoFix_n(f_c))}{n}$$ exists, and moreover is contained in the interval  $[\log 2.0890, \log 2.1071]$. That the polynomial $2t^3-3(t^2-1)-4t$ in part 2 of Conjecture \ref{Conjecture3} is irreducible over $\mathbb{Q}$ can be easily checked by many means. (For example, if it were to be reducible over $\mathbb{Q}$, then it would have at least one root in $\mathbb{Q}$, which then must be half of an integer. Using computer softwares such as Mathematica, one finds that the given polynomial has 3 real roots lying between $-2.5$ and $2.5$, and no half integer in this interval is a root of the polynomial.) 

\section{Method} The method in this paper is a mixture of rigorously theoretical and experimental ones. 

We use the current techniques in complex Dynamical Systems to prove the mentioned relations between the conjectures in the Introduction section. We also use rigorously theoretical arguments to show that, provided the pattern observed in our experiments hold for all iterate, then our maps (on the smooth affine quadric 3-folds) are primitive, and to reduce the computation of dynamical degrees on the quadric 3-folds to that on the affine $4$-space.  

On the other hand, the current techniques in complex Dynamical Systems are not enough for to resolve  Conjecture \ref{Conjecture3}, see the Discussion section and the Conclusion section for more detail. Hence, we need to rely on computational methods. For systems of polynomial equations over the field of complex numbers, there are symbolic methods such as Gr\"obner basis which have strong theoretical guarantees but usually cannot find individual solutions to the system (which is what we need here since we want to check whether a periodic point is hyperbolic or not) and usually is too slow, in particular for the systems of equations related to computing periodic points of our map. Hence, we need to utilise numerical methods, more specifically here the software Bertini \cite{bates-etal}. 

Here is a brief summary of the idea for using experimental results reported later in this paper to give support to the validity of Conjecture \ref{Conjecture3}. {To check} part 2 of Conjecture \ref{Conjecture3}, we compare the degree sequences of the iterates of $F$ and some sequences related to the linear recurrence $w_n=3/2 \cdot (w_{n-1}-w_{n-3})+2w_{n-2}$, whose exponential growth is $\zeta _2$. For to check the Weaker estimate in part 3 of Conjecture \ref{Conjecture3}, we test if the sequence $  (\sharp IsoFix_{2n+1}(f_c))^{1/(2n+1)}$ is {\bf decreasing}, which is related to the log concavity phenomenon in Proposition \ref{TheoremMain} (see next), and to observe that for $n=4,5$ the number $  (\sharp IsoFix_{2n+1}(f_c))^{1/(2n+1)}$ is $\leq 2.108$. For to check the Stronger estimate in part 3 of Conjecture \ref{Conjecture3}, we compare the two sequences $(\sharp IsoFix_{2n+1}(f_c))^{1/(2n+1)}$ and $(\sharp IsoFix_{2n}(f_c))^{1/(2n)}$ and use the Weaker estimate.  

The next result provides a heuristic support for the mentioned decreasing phenomenon of the sequence $(\sharp IsoFix_{2n+1}(f_c))^{1/(2n+1)}$. To see why, we recall that roughly speaking (under the assumption that $g$ has only isolated fixed points)  the Lefschetz number $L(g)$ is the sum of fixed points (with multiplicities) of $g$.  Hence, the observed decreasing phenomenon would be a consequence of the next result, if it were the case that all fixed points of $\widehat{f_c}^n$ were in $Z_c$ and isolated. 

\begin{proposition} Let $X\subset \mathbb{P}^{2m}$ be a smooth quadric. Let $g:X\dashrightarrow X$ be a dominant rational map. Then the Lefschetz numbers $L(g^n)$ (i.e. the intersection number between the graph of $g^n$ and the diagonal of $X$) are all positive, and satisfy log concavity, i.e.
\begin{eqnarray*}
L(g^{n+n'})\leq L(g^n)L(g^{n'})
\end{eqnarray*}
for all $n,n'\geq 0$. In particular, the sequence $\{L(g^n)^{1/n}\}_{n=1,2,\ldots}$ is decreasing. 

Similarly, for every $0\leq j\leq 2m-1$, the degree sequence $\{||(g^n)^*|_{H^{2j}(X)}||\}_{n=1,2,\ldots }$ (where $||.||$ is a given norm on $H^{2j}(X)$) satisfies the log concavity. 

\label{TheoremMain}\end{proposition}

\section{Results}

In this section we present our results around the map $f$ and Conjecture \ref{Conjecture3}. These results consist of both rigorously theoretically proven ones and experimental ones. Based on these results, we propose a road map toward resolving Conjecture~\ref{Conjecture3} in the affirmative. 

Let $F:~\mathbb{C}^4\rightarrow \mathbb{C}^4$, $X_c$, $Z_c$ and $f_c:Z_c\rightarrow Z_c$ be defined as in the introductive section.  Throughout this section, we use the following notations: 

$\widehat{F}=$ the extension of $F$ to the projective space $\mathbb{P}^4$;

$\widehat{f}_c=$ the extension of $f_c$ to $X_c$;

$d_n^{(1)}=$ the first degree of $F^{n}$ (the n-th iterate of $F$): it is defined as the degree of the inverse image by $\widehat{F}^n$ of  a generic linear 3-dimensional subspace $H$ in $\mathbb{P}^4$; 

$d_n^{(2)}=$ the second degree of $F^{n}$: it is defined as the degree of the inverse image by $\widehat{F}^n$ of  a generic linear 2-dimensional subspace $H^2$ in $\mathbb{P}^4$;

$d_n^{(3)}=$ the third degree of $F^{n}$: it is defined as the degree of the inverse image by $\widehat{F}^n$ of  a generic linear 1-dimensional subspace $H^3$ in $\mathbb{P}^4$ (it is indeed the same as the first degree of the inverse map $F^{-n}$);

$b_n$ is the sequence satsifying the linear recurrence $b_n=3/2 \cdot (b_{n-1}-b_{n-3})+2b_{n-2}$, with the initial terms $b_1=3$, $b_2=7$ and $b_3=17$;

$c_n$ is the sequence satsifying the linear recurrence $c_n=3/2 \cdot (c_{n-1}-c_{n-3})+2c_{n-2}$, with the initial terms $c_1=5$, $c_2=9$ and $c_3=25$;

$IsoFix_n=$ the set of isolated fixed points of $f_c^n$;

$C\subset Z_c$: the curve defined by the ideal $<x_2-x_1^2x_2-x_3,x_1+x_4, x_1x_4-x_2x_3-c>$. 

$D_1\subset Z_c$: the curve with 2 components defined by the ideals $<x_2-x_1^2x_2-x_3,x_1+x_4, x_1x_4-x_2x_3-c>$ and $<-x_2+x_1^2x_2-x_3,x_1-x_4, x_1x_4-x_2x_3-c>$;

$D_2\subset Z_c$: the curve with 2 components defined by the ideals  $<x_2, x_3, x_1x_4-x_2x_3-c>$ and $<x_1,x_4, x_1x_4-x_2x_3-c>$. 
 
Here are some useful comments concerning the above notations. First, by the definition of dynamical degrees, we have 
\begin{eqnarray*}
\lim _{n\rightarrow\infty}[d_n^{(1)}]^{1/n}&=&\lambda _1(\widehat{F}),\\
\lim _{n\rightarrow\infty}[d_n^{(2)}]^{1/n}&=&\lambda _2(\widehat{F}),\\
\lim _{n\rightarrow\infty}[d_n^{(3)}]^{1/n}&=&\lambda _3(\widehat{F}).
\end{eqnarray*}

Since $\widehat{F}$ is birational,  we have $\lambda _0(\widehat{F})=\lambda _4(\widehat{F})=1$. Similarly, we have $\lambda _0(f_c)=\lambda _3(f_c)=1$. The relation between the dynamical degrees of $\widehat{F}$ and those of $\widehat{f_c}$ will be given in the next Subsection. 

From our experiments, to be presented later in this section, the fixed point set of $f_c^{n}$ seems to have the following structure: 

The fixed point set of $f_c^{4n+2}$ is the union of the curve C and $IsoFix_{4n+2}(f_c)$;

The fixed point set of $f_c^{4n}$ is the union of the curve D (which is the union of $D_1$ and $D_2$ - of different multiplicities) and $IsoFix_{4n}(f_c)$;

The fixed point set of $f_c^{2n+1}$ is $IsoFix_{2n+1}(f_c)$; 

All isolated periodic points of $f_c$ are hyperbolic, i.e. $IsoFix_{n}(f_c)=HP_n(f_c)$.

\subsection{Theoretical results}

We first start with some relations between Conjectures \ref{Conjecture1}, \ref{Conjecture2}, \ref{ConjectureAlgebraicInteger} and  \ref{Conjecture3}. 

\begin{lemma}
1) If Conjecture \ref{Conjecture1} holds, then Conjecture \ref{Conjecture2} holds. 

2) Assume that parts 1 and 2 of Conjecture \ref{Conjecture3} hold. 

a) If the Stronger estimate of part 3 of Conjecture \ref{Conjecture3} holds, then Conjecture \ref{Conjecture2} does not hold. 

b) If the Weaker estimate of part 3 of Conjecture \ref{Conjecture3} holds, then Conjecture \ref{Conjecture1} does not hold.  

3) If part 2 of Conjecture \ref{Conjecture3} holds, then Conjecture \ref{ConjectureAlgebraicInteger} does not hold. 

\end{lemma}
\begin{proof}

1) Assume that Conjecture \ref{Conjecture1} holds. Let $f:X\dashrightarrow X$ be a dominant rational map which is cohomologically hyperbolic. Let $\pi :Y\dashrightarrow X$ be the birational map given by Conjecture \ref{Conjecture1}. 

Let $Z\subset X$ be a Zariski open dense set.  Since $\pi :Y\dashrightarrow X$ is birational, there is a Zariski open dense set $U\subset Z$ and a Zariski open dense set $V\subset Y$ so that $\pi$ induces an isomorphism between $V$ and $U$. 

By Conjecture \ref{Conjecture1}, $HP_n(Y)$ equidistributes to $\mu$ and where $\mu$ has no mass on proper analytic subsets. Since $Y\backslash V$ is a proper analytic subset of $Y$, it follows that $HP_n(f_Y)\cap V$ equidistributes to $\mu$ as well, and the exponential growth of $\sharp [HP_n(f_Y)\cap V]$ is also $\lambda (f_Y)$. Since $\pi :V\rightarrow U$ is an isomorphism, it follows that $HP_n(f_Y)\cap V=HP_n(f)\cap U$. Since $\lambda (f_Y)=\lambda (f)$ by the birational invariance of dynamical degrees, we obtain: 
\begin{eqnarray*}
\lim _{n\rightarrow\infty}\frac{\log \sharp [HP_n(f)\cap U]}{n}=\log \lambda (f). 
\end{eqnarray*}

Since $HP_n(f)\cap U\subset IsoFix_n(f)\cap Z$, we get: 
\begin{eqnarray*}
\liminf _{n\rightarrow\infty}\frac{\log \sharp  [IsoFix_n(f)\cap Z]}{n}\geq \log \lambda (f). 
\end{eqnarray*}

On the other hand, by \cite{dinh-etal2}, we have 
\begin{eqnarray*}
\limsup _{n\rightarrow\infty}\frac{\log \sharp  IsoFix_n(f)}{n}\leq  \log \lambda (f).
\end{eqnarray*}

Combining the above two inequalities, we have finally
\begin{eqnarray*}
\lim _{n\rightarrow\infty}\frac{\log \sharp  [IsoFix_n(f)\cap Z]}{n}=\log \lambda (f).
\end{eqnarray*}
This means that Conjecture \ref{Conjecture2} holds. 

2) We prove for part a) only, the proof of part b) is similar. Assume that parts 1 and 2 of Conjecture \ref{Conjecture3} hold, and also that the Stronger estimate of part 3 of Conjecture \ref{Conjecture3} holds. Since $X_c$ is of dimension $3$, and the dynamical degrees of $\widehat{f}_c$ are: $\lambda _0(\widehat{f}_c)=1$, $\lambda _1(\widehat{f}_c)=\zeta _1\sim 1.8393$, $\lambda _2(\widehat{f}_c)=\zeta _2\sim 2.1108$, and $\lambda _3(\widehat{f}_c)=1$, the map $\widehat{f}_c$ is cohomologically hyperbolic. We also have that $Z_c$ is a Zariski open dense set of $X_c$, and $\widehat{f}_c|_{Z_c}=f_c$. Hence, $IsoFix_n(\widehat{f}_c)\cap Z_c=IsoFix_n(f_c)$. By Conjecture \ref{Conjecture3} we have 
\begin{eqnarray*}
\limsup _{n\rightarrow\infty}\frac{\log \sharp IsoFix_n(f_c)}{n}\leq \log 2.108 < \log  \zeta _2=\log \lambda (f_c). 
\end{eqnarray*}
This contradicts Conjecture \ref{Conjecture2}. 

3) Assume that part 2 of Conjecture \ref{Conjecture3} holds. Then we have $\lambda _2(\widehat{f}_c)=\zeta _2$ is an algebraic number, but not an algebraic integer. By Lemma \ref{Lemma2}, we have $\lambda _3(\widehat{F})=\lambda _2(\widehat{f}_c)$, hence $\lambda _3(\widehat{F})$ is not an algebraic integer. In this case, also $\lambda _2(\widehat{F})=\zeta _2$ is not an algebraic integer. Thus the polynomial map $F:\mathbb{C}^4\rightarrow \mathbb{C}^4$ is a counter-example to Conjecture \ref{ConjectureAlgebraicInteger}. 
\end{proof}

We recall that \cite{zhang} a dominant rational map $f:X\dashrightarrow X$ is primitive, if there {\bf do not exist} a variety $W$ of dimension $1\leq \dim (W) \leq \dim (X)-1$, a dominant rational map $\pi : X\dashrightarrow W$, and a dominant rational map $g:W\dashrightarrow W$ so that $\pi \circ f=g\circ \pi$. Primitive maps can be viewed as building blocks from which all maps can be constructed. It is clear from the definition that being primitive is a birational invariant. We have the following result. 

\begin{lemma}
Assume that parts 1 and 2 of Conjecture \ref{Conjecture3} hold. Then the map $f_c$ is primitive. 
\label{Lemma3}\end{lemma}
\begin{proof}
Let $\widehat{f}_c^{-1}$ be the inverse of $\widehat{f}_c$. If Conjecture \ref{Conjecture3} holds, then 
\begin{eqnarray*}
\lambda _1(\widehat{f}_c^{-1})=\lambda _2(\widehat{f}_c)>\lambda _1(\widehat{f}_c)=\lambda _2(\widehat{f}_c^{-1}).
\end{eqnarray*}
This inequality implies that $\widehat{f}_c^{-1}$ is primitive, \cite{oguiso-truong}. Hence $f_c^{-1}$ (the inverse of $f_c$) is also primitive.

From this, we will show that $f_c$ is primitive. Assume by contradiction that $f_c$ is not primitive. Then there are $\pi :X\dashrightarrow W$ and $g:W\dashrightarrow W$ dominant rational maps (with $1\leq \dim (W)\leq 2$) so that $\pi \circ f_c=g\circ \pi$. Since $f_c$ is a birational map, it is easy to see that $g$ is also a birational map. Then, it follows that $\pi \circ f_c^{-1}=g^{-1}\circ \pi$, which contradicts the fact that $f_c^{-1}$ is primitive.  
\end{proof}

{We next relate} dynamical degrees of $\widehat{F}$ and those of $\widehat{f}_c$. 

\begin{lemma}
We have

1) $\lambda _1(\widehat{f}_c)=\lambda _1(\widehat{F})$, and $\lambda _2(\widehat{f}_c)=\lambda _3(\widehat{F})$. 

2) $\lambda _2(\widehat{F})=\max \{\lambda _1(\widehat{F}), \lambda _3(\widehat{F})\}$.

\label{Lemma2}\end{lemma}
\begin{proof}

We know from the introductive section that $F$ is semi-conjugated to the identity map $id_{\mathbb{C}}$ on the curve $\mathbb{C}$, via the map $\phi : \mathbb{C}^4\rightarrow \mathbb{C}$ with $\phi (x_1,x_2,x_3,x_4)=x_1x_4-x_2x_3$. The dynamical degrees of (the extension to $\mathbb{P}^1$ of)  $id_{\mathbb{C}}$ are $\lambda _0(id_{\mathbb{C}})=\lambda _1(id_{\mathbb{C}})=1$.   

Let $\lambda _j(\widehat{F}|\phi )$ ($j=0,1,2,3$) be the relative dynamical degrees w.r.t. $\phi$ of $\widehat{F}$ \cite{dinh-nguyen} \cite{dinh-etal3}. Since $F$ preserves the fibres of $\phi$, it follows that $\lambda _j(\widehat{F}|\phi  )=\lambda _j(\widehat{f}_c)$.  

By \cite{dinh-nguyen} \cite{dinh-etal3}, we have 
\begin{eqnarray*}
\lambda _1(\widehat{F})&=&\max \{\lambda _0(id _{\mathbb{C}})\lambda _1(\widehat{f}_c),\lambda _1(id _{\mathbb{C}})\lambda _0(\widehat{f}_c)\}=\lambda _1(\widehat{f}_c),\\
\lambda _3(\widehat{F})&=&\max \{\lambda _0(id _{\mathbb{C}})\lambda _3(\widehat{f}_c),\lambda _1(id _{\mathbb{C}})\lambda _2(\widehat{f}_c)\}=\lambda _2(\widehat{f}_c),\\
\lambda _2(\widehat{F})&=&\max \{\lambda _0(id _{\mathbb{C}})\lambda _2(\widehat{f}_c),\lambda _1(id _{\mathbb{C}})\lambda _1(\widehat{f}_c)\}\\
&=&\max\{\lambda _1(\widehat{f}_c), \lambda _2(\widehat{f}_c)\}=\max\{\lambda _1(\widehat{F}), \lambda _3(\widehat{F})\}. 
\end{eqnarray*}

\end{proof}

Next, we compute $\lambda _1(\widehat{f_c})$. 

\begin{lemma}
$\lambda _1(\widehat{f}_c)=\zeta _1$, the largest root of the polynomial $t^3-t^2-t-1$.  
\label{Lemma4}\end{lemma}
\begin{proof} By Lemma \ref{Lemma2}, we have $\lambda _1(\widehat{f}_c)=\lambda _1(\widehat{F})$.  Thus we only need to show that $\lambda _1(\widehat{F})=\zeta _1$. Recall that $\lambda _1(\widehat{F})=\lim _{n\rightarrow\infty}[d_n^{(1)}]^{1/n}$. Hence, the proof is finished if we can show that the degree sequence $d_n^{(1)}$ of $F$ satisfies the linear recurrence: 
\begin{eqnarray*}
d_{n}^{(1)}=d_{n-1}^{(1)}+d_{n-2}^{(1)}+d_{n-3}^{(1)},
\end{eqnarray*}
for all $n$. 

The leading term of $F(x_1,x_2,x_3,x_4)=((F)_1,(F)_2,(F)_3,(F)_4)$ in terms of degree is $x_1x_2x_4$ in $(F)_4.$ Moreover, $(F)_3$ has a unique leading monomial, and $x_1x_2x_4$ is the unique leading monomial in $(F)_4$. By direct calculation, we find that $F^2=((F^2)_1,(F^2)_2,(F^2)_3,(F^2)_4)$ is given by: 
$$(-x_4,x_3-x_1x_2x_4,x_2-x_2x_4^2,-x_1+x_1x_2^2+x_2x_3x_4-x_1x_2^2x_4^2).$$
Hence, again, the leading term of $F^2$ is contained in $(F^2)_4$, and $(F^2)_4$ has a unique leading monomial. One can directly check the same phenomenon for $F^3$ and $F^4$. 

We will prove by induction that the leading term of $F^n=((F^n)_1, (F^n)_2, (F^n)_3,(F^n)_4)$ is contained in $(F^n)_4$, and moreover $(F^n)_4$ has a unique leading monomial. In addition the concerned linear recurrence for the degree sequence holds. Assume by induction that {this claim holds} until $n-1$. We will show that it also holds for $n$. Since $F^{n}(x_1,x_2,x_3,x_4)$ is
$$((F^{n-1})_2, -(F^{n-1})_4, (F^{n-1})_1 - (F^{n-1})_1
\cdot{(F^{n-1})_2}^2, -(F^{n-1})_3+(F^{n-1})_1\cdot(F^{n-1})_2 \cdot (F^{n-1})_4)$$
and since 
$$(F^{n-1})_2 =-(F^{n-2})_4$$ and 
$$(F^{n-1})_1 = (F^{n-2})_2=-(F^{n-3})_4$$ it follows that the leading degree of $F^n$ is that of $(F^{n-1})_4\cdot(F^{n-2})_4 \cdot (F^{n-3})_4$ which is contained in $(F^n)_4$, and that $(F^n)_4$ has a unique leading monomial, which is the product of the unique leading monomials in $(F^{n-1})_4$, $(F^{n-2})_4$ and $(F^{n-3})_4$. From this,  we have immediately the linear recurrence: 
\begin{eqnarray*}
d_{n}^{(1)}=d_{n-1}^{(1)}+d_{n-2}^{(1)}+d_{n-3}^{(1)}.
\end{eqnarray*}
\end{proof}

Finally, we prove Proposition \ref{TheoremMain}. 
\begin{proof}[Proof of Proposition \ref{TheoremMain}]

First, we have the known fact that the cohomology groups of $X$ come from the pullback on cohomology of the embedding $\iota :X\hookrightarrow \mathbb{P}^{2m}$. [For the convenience of the readers, we briefly recall the arguments. By the Lefschetz hyperplane theorem, for $j\leq 2m-2$, the pullback $\iota ^*: H^{j}(\mathbb{P}^{2m})\rightarrow H^{j}(X)$ is an isomorphism. For $j\geq 2m+2$ one can use Poincare duality to determine the cohomology group. It remains to show that $H^{2m-1}(X)=0$. This can be done by computing the Euler characteristic of $X$, and the latter can be done by using the normal bundle sequence and Whitney sum formula.]

In particular, this means the following: $H^{2j+1}(X)=0$ for all $j$; if $2j$ is even and $\leq 2m-2$ then $H^{2j}(X)$ is generated by $h^j$ where $h=\iota ^*(H)$ is the pullback of a hyperplane $H$ in $\mathbb{P}^{2m}$, while the remaining groups are computed using Poincare duality. Hence $H^*(X)$ is generated by $h$, and we have $h^{2m-1}=H^{2m-1}.X=2$. 

Thus, if $\pi _1,\pi _2:X\times X\rightarrow X$ are the two canonical projections, then the cohomology class of the diagonal of $X$ is 
\begin{eqnarray*}
\{\Delta _X\}=\frac{1}{2}\sum _{j=0}^{2m-1}\pi _1^*(h^j).\pi _2^*(h^{2m-1-j}).  
\end{eqnarray*}
{Indeed, by Kunneth's formula and the description of the cohomology group of $X$ from the previous paragraph, we have $\{\Delta _X\}=\sum _{j=0}^{2m-1}a_j\pi _1^*(h^j).\pi _2^*(h^{2m-1-j})$, for constants $a_j$. To determine $a_i$, we can do as follows. By the definition of the diagonal, we have $(\pi _1)_*(\{\Delta _X\}.\pi _2^*(h^i))=h^i$ for all $0\leq i\leq 2m-1$. On the other hand, from the formula for $\{\Delta _X\}$, by the projection formula, we obtain easily $(\pi _1)_*(\{\Delta _X\}.\pi _2^*(h^i))=(\pi _1)_*(a_i\pi _1^*(h^i).\pi _2^*(h^{2m-1-i}).\pi _2^*(h^{i}))=a_i(h^{2m-1}).h^i$. Since $h^{2m-1}$ is $2$, it follows that $2a_i=1$ for all $i$, and hence $a_i=1/2$ for all $i$.} 

From this, we get the following formula for the Lefschetz number of a map $g:X\rightarrow X$: 
\begin{eqnarray*}
L(g)=\sum _{j=0}^{2m-1}g^*\left(\frac{1}{2}h^j\right).h^{2m-1-j}.
\end{eqnarray*}
{In fact, by the previous formula for $\{\Delta_X\}$, we have 
$$L(g)=\{\Gamma _g\}.\{\Delta _X\} =\sum _{j=0}^{2m-1}\frac{1}{2}\{\Gamma _j\}.\pi _1^*(h^j).\pi _2^*(h^{2m-1-j}).$$ 
For each $j$, we have 
$$\frac{1}{2}\{\Gamma _j\}.\pi _1^*(h^j).\pi _2^*(h^{2m-1-j})= \frac{1}{2}(\pi _2)_*(\{\Gamma _j\}.\pi _1^*(h^j)).h^{2m-1-j}.$$ 
By definition, we have $(\pi _2)_*(\{\Gamma _j\}.\pi _1^*(h^j))=g^*(h^j)$, thus we obtain the claimed expression for $L(g)$.}
{Each of the summands} in the right hand side is $>0$, and hence $L(g)>0$. We denote by $d_{j}(g):=g^*(\frac{1}{2}h^j).h^{2m-1-j}$ the $j$-th degree of $g$. 

Now, we will show that for any two dominant rational maps $g_1,g_2:X\dashrightarrow X$, and any $j\leq 2m-1$ then $d_j(g_1\circ g_2)\leq d_j(g_1)d_j(g_2)$. To this end, we follow the ideas in \cite{russakovskii-shiffman} \cite{dinh-sibony10} of using automorphisms of $X$ to regularise positive closed currents. We recall that the quadric $X$ is a homogeneous space, whose automorphism group is the subgroup of the linear automorphisms of $\mathbb{P}^{2m}$ preserving the quadratic form. Hence, by using a  convolution process with the aid of the Haar measure of the automorphism group of $X$, for any positive closed $(j,j)$ current $T$ on $X$, there is a sequence of positive closed smooth $(j,j)$ forms $T_{\epsilon}$ weakly converging to $T$ so that (recall that the cohomology of $X$ is generated by $h$) in cohomology $\{T_{\epsilon}\}=\{T\}$. 

Now we recall how $d_j(g_1\circ g_2)$ can be computed. There is a Zariski open set $Z\subset X$ so that if we choose  generic algebraic varieties in $X$: $V_j$ of codimension $j$ representing the cohomology class $h^j$ and $W_{2m-j}$ of codimension $2m-j$ representing the cohomology class $h^{2m-j}$, then 
\begin{eqnarray*}
d(g_1\circ g_2)=(g_2|_{Z})^{-1}(g_1^{-1}(V_j/2)).W_{2m-j}.
\end{eqnarray*}
Here $g_1|_Z$ is a proper map of finite fibres, and hence the preimage of $g_1|_Z$ of any variety is again a variety of the same dimension. Moreover, we have the following property: If $T_{\epsilon}$ is a sequence of positive closed smooth forms weakly converging to the current of integration over $g_1^{-1}(V_j)/2$ and being of the same cohomology class as that of $g_1^{-1}(V_j)/2$ (which is proven in the above), then
\begin{eqnarray*}
(g_2|_{Z})^{-1}(g_1^{-1}(V_j/2)).W_{2m-j}\leq \lim _{\epsilon\rightarrow 0} g_2^*(\{T_{\epsilon }\}).h^{2m-j}, 
\end{eqnarray*}
 {because $g_2^*(\{T_{\epsilon }\}).h^{2m-j}$ is the cohomology class of  $g_2^*(T_{\epsilon}).W_{2m-j}$, and the latter is a positive closed current for all $\epsilon$, and any cluster point of $\{g_2^*(T_{\epsilon}).W_{2m-j}\}$ (when $\epsilon$ decreases to $0$) will coincide with $(g_2|_{Z})^{-1}(g_1^{-1}(V_j/2)).W_{2m-j}$ when restricted to $Z$}. 
Since 
\begin{eqnarray*}
\{T_{\epsilon }\}&=&\{g_1^*(V_j/2)\}=g_1^*(h^j/2)\\
&=&[g_1^*(h^j/2).h^{2m-j}]/2.h^{j}=d_j(g_1)\frac{h^j}{2}, 
\end{eqnarray*}
we obtain 
\begin{eqnarray*}
(g_2|_{Z})^{-1}(g_1^{-1}(V_j/2)).W_{2m-j}\leq d_j(g_1)g_2^*(h^j/2).h^{2m-j}=d_j(g_1)d_j(g_2). 
\end{eqnarray*}
Therefore, the degree sequences are log concave. From this, we obtain
\begin{eqnarray*}
L(g_1\circ g_2)&=&\sum _jd_j(g_1\circ g_2)\leq \sum _jd_j(g_1)d_j(g_2)\\
&\leq&(\sum _jd_j(g_1))(\sum _jd_j(g_2))=L(g_1)L(g_2). 
\end{eqnarray*}

Applying for $g_1=g^{n}$ and $g_2=g^{n'}$ we obtain the log concavity $L(g^{n+n'})\leq L(g^n)L(g^{n'})$ needed. From this log concavity property, it is well known that the sequence $n\mapsto L(g^n)^{1/n}$ is decreasing. 
\end{proof}

\begin{remark}
One can also prove Proposition \ref{TheoremMain} by a purely algebraic proof, which replaces regularisation of currents by Chow's moving lemma and which is valid on any algebraically closed field, see \cite{truong}.
\end{remark}

\subsection{Experimental results} In this subsection we calculate the fixed point sets of the iterates $f_c^n$, as well as the degree sequences $d_n^{(1)}$, $d_n^{(2)}$ and $d_n^{(3)}$, for as large as possible $n$'s, to help study Conjecture \ref{Conjecture3}. We also study how close the degree sequence $d_n^{(2)}$ and $d_n^{(3)}$ are to the linear recurrence $v_n=3/2 \cdot (v_{n-1}-v_{n-3})+2v_{n-2}$ related to the polynomial $t^3-3/2 \cdot (t^2-1)-2t$ (for which $\zeta _2$ is the largest root). These calculations will be used in the next subsection, where we propose an approach towards solving Conjecture \ref{Conjecture3}. Here we recall the relevant relations from Lemmas \ref{Lemma2} and \ref{Lemma4}: 
\begin{eqnarray*}
\lim _{n\rightarrow\infty}[d_n^{(1)}]^{1/n}&=&\zeta _1=\lambda _1(\widehat{F})=\lambda _1(\widehat{f}_c),\\
\lim _{n\rightarrow\infty}[d_n^{(3)}]^{1/n}&=&\lambda _3(\widehat{F})=\lambda _2(\widehat{f}_c),\\
\lim _{n\rightarrow\infty}[d_n^{(2)}]^{1/n}&=&\lambda _2(\widehat{F})=\max \{\lambda _1(\widehat{F}),\lambda _3(\widehat{F})\}=\max \{\lambda _1(\widehat{f}_c),\lambda _2(\widehat{f}_c)\}.
\end{eqnarray*}

Formal computer algebra techniques (such as Gr\"obner basis routine on the softwares Mathematica and Maple) can only compute up to about $N=5$. Hence, we
utilize regeneration~\cite{Regen}\cite{RegenCascade}\cite{RegenExtension} and the trace test~\cite{Trace}\cite{Trace2}\cite{Trace1}
implemented in Bertini~\cite{bates-etal}.
{Although numerical routines, the
trace test provides a high level of confidence 
with reliable performance.  
Moreover, parallel computing
can be used to speed up the computations}.
 We note that the computations for finding the periodic points are more expensive and difficult than those for calculating the degree sequences.  

The first table reports on the calculations for periodic points of $f_c$ (up to period $n=12$), as well as the exponential growth of isolated periodic points. All the isolated periodic points are hyperbolic. 

$$\begin{array}{c|c|c}
N & \hbox{fixed points on general fiber}&[\sharp IsoFix_{N}(f_c)]^{1/N}\\ \hline 
1 & 4&{\bf 4} \\
2 & C ~(\hbox{occurring with multiplicity 1})&0\\
3 & 10&{\bf 2.15443469003} \\
4 & D_1~ (\hbox{multiplicity 1})~ $\&$ ~D_2 ~(\hbox{multiplicity 2}) &0\\
5 & 44&{\bf 2.13152551327} \\
6 & C ~(\hbox{multiplicity 1}) \hbox{~AND~ 12 points}&1.51308574942 \\
7 & 186&{\bf 2.10967780991} \\
8 & D_1~ (\hbox{multiplicity 1})~ $\&$ ~D_2 ~(\hbox{multiplicity 2}) \hbox{~AND~128 points}&1.83400808641\\
9 & 820&{\bf 2.10744910267} \\
10 & C ~(\hbox{multiplicity 1}) \hbox{~AND~ 1440 points}&2.06936094886 \\
11 & 3634& {\bf 2.10703309279} \\
12 & D_1~ (\hbox{multiplicity 1})~ $\&$ ~D_2 ~(\hbox{multiplicity 2}) \hbox{~AND~6908 points}&2.08903649661\\
\end{array}$$

The next table computes the degree sequence for the iterates $F^n$ (up to $n=14$), as well as their exponential growth: 

$$\begin{array}{c|c|c|c|c|c|c}
N & d_N^{(1)} & d_N^{(2)} & d_N^{(3)}&[d_N^{(1)}]^{1/N}&[d_N^{(2)}]^{1/N}&[d_N^{(3)}]^{1/N}\\ \hline
1 & 3 & 5 & 3&3&5&3 \\
2 & 5 & 9 & 7&2.2360679775&3&2.64575131106 \\
3 & 9 & 25 & 17&2.08008382305&2.92401773821&2.57128159066\\
4 & 17 & 49 & 37&2.03054318487& 2.64575131106& 2.46632571456\\
5 & 31 & 109 & 79&1.98734075466& 2.55555539674&2.39621299048\\
6 & 57 & 225 & 167&1.96175970274&2.46621207433& 2.34667391139\\
7 & 105 & 477 & 353&1.94420174432&2.41348988334&2.3118934527 \\
8 & 193 & 1005 & 745&1.93061049898& 2.37285258221&2.28570160944\\
9 & 355 & 2117 & 1571&1.92025412137&2.34166378698&2.26532588341\\
10 & 653 & 4465 & 3311&1.91201510161&2.31729938473& 2.24903346712\\
11 & 1201 & 9401 & 6977&1.90527844956&2.29719383004& 2.23575612581\\
12 & 2209 & 19817 & 14701&1.8996910486&2.28079626154&2.22473817189 \\
13 & 4063 & 41741 & 30975&1.89497551023&2.2668767672&2.2154523255 \\
14 & 7473 & 87961 & 65263&1.89094202127&2.25508846088&2.2075207175 \\
\end{array}$$

We can see from the above two tables that, as predicted by Conjecture \ref{Conjecture3}, the exponential growth of the degree sequence is much larger than that of the isolated periodic points. From the table for the degree sequence above, we can readily check that the sequence $d_n^{(1)}$'s indeed satisfies the linear recurrence $w_n=w_{n-1}+w_{n-2}+w_{n-3}$, as proven in Lemma \ref{Lemma4}. If Conjecture \ref{Conjecture3} holds, then we must have $\lambda _2(\widehat{F})=\lambda _3(\widehat{F})$, and a first idea towards actually showing that $\lambda _2(\widehat{F})=\lambda _2(\widehat{F})=\zeta _2$ is to show that both sequences $d_n^{(2)}$'s and $d_n^{(3)}$'s satisfy the linear recurrence $w_n=3/2 \cdot (w_{n-1}-w_{n-3})+2w_{n-2}$. It turns out that neither of these sequences satisfies this linear recurrence, but the next two tables show that the differences to this linear recurrence are  relatively small (in comparison to the size of the concerned degree sequences).

 $$\begin{array}{c|c|c|c}
N & d_N^{(2)}& 3/2 \cdot (d_{N-1}^{(2)}-d_{N-3}^{(2)})+2d_{N-2}^{(2)} &d_N^{(2)}-[3/2 \cdot (d_{N-1}^{(2)}-d_{N-3}^{(2)})+2d_{N-2}^{(2)}] \\ \hline
1 & 5 & & \\
2 & 9 & & \\
3 & 25 & & \\
4 & 49 & 48 & 1 \\ 
5 & 109 & 110 & -1\\
6 & 225 & 224 &1 \\
7 & 477 & 482 & -5 \\
8 & 1005 &1002 & 3\\ 
9 & 2117 & 2124 & -7\\
10 &4465  & 4470 & -5\\
11 &9401  & 9424 & -23\\
12 & 19817 & 19856 &-39 \\ 
13 &  41741&41830  &-89 \\
14 & 87961 & 88144 & -183
\end{array}$$

    $$\begin{array}{c|c|c|c}
N & d_N^{(3)}& 3/2 \cdot (d_{N-1}^{(3)}-d_{N-3}^{(3)})+2d_{N-2}^{(3)} & d_N^{(3)}-[3/2 \cdot (d_{N-1}^{(3)}-d_{N-3}^{(3)})+2d_{N-2}^{(3)} ] \\ \hline
1 & 3 & & \\
2 & 7 & & \\
3 & 17 & & \\
4 & 37 & 35 & 2 \\ 
5 & 79 & 79 & 0\\
6 & 167 & 167 & 0\\
7 & 353 & 353 & 0\\
8 & 745 & 745 & 0\\ 
9 & 1571 & 1573 & -2 \\
10 & 3311 & 3317 & -6 \\
11 & 6977 & 6991 & -14 \\
12 & 14701 & 14731 & -30\\ 
13 & 30975 & 31039 & -64 \\
14 & 65263 & 65399 & -136 
\end{array}$$

\subsection{Further analysis $\&$ A  road map towards an affirmative answer to Conjecture \ref{Conjecture3}}

In this Subsection we further analyze the experimental findings in the previous Subsection, in connection to Conjecture \ref{Conjecture3}. Since part 1 of Conjecture \ref{Conjecture3} is solved by Lemma \ref{Lemma4}, it remains to treat parts 2 and 3 of Conjecture \ref{Conjecture3}. 

\subsubsection{An approach towards establishing the upper bound $\lambda _2(\widehat{f}_c)\leq \zeta _2$} 
As seen from before, this is equivalent to establishing that 
\begin{eqnarray*}
\lim _{n\rightarrow\infty}[d_n^{(3)}]^{1/n}\leq \zeta _2.
\end{eqnarray*}

Here is an approach to showing this. From the last 2 tables in the previous Subsection, it seems very {plausible} that for $n\geq 5$, we should have $d_n^{(3)}\leq 3/2 \cdot (d_{n-1}^{(3)}-d_{n-3}^{(3)})+2d_{n-2}^{(3)}$.  Also, it seems very evident that $d_{n+1}^{(3)}\geq 2d_n^{(3)}$ for all $n\geq 1$. This is indeed enough to proving the desired upper bound, as seen in the next lemma.

\begin{lemma}
Assume that for all $n\geq 5$, we have
\begin{eqnarray*}
d_n^{(3)}&\leq& 3/2 \cdot (d_{n-1}^{(3)}-d_{n-3}^{(3)})+2d_{n-2}^{(3)},\\
d_{n}^{(3)}&\geq & 2 d_{n-1}^{(3)}. 
\end{eqnarray*}
Then $\lambda _3(\widehat{F})\leq \zeta _2$. 
\label{Lemma5}\end{lemma}
\begin{proof} Because $[d_n^{(3)}]^{1/n}$ is a decreasing sequence (see  \cite{russakovskii-shiffman}, and see also Proposition \ref{TheoremMain}), we have $d_n^{(3)}/d_{n-1}^{(3)}\leq [d_{n-1}^{(3)}]^{1/(n-1)}$ for all $n\geq 2$.  Therefore, from the assumption that $d_{n}^{(3)}\geq  2 d_{n-1}^{(3)}$ for all $n\geq 5$, we obtain
\begin{eqnarray*}
2\leq \limsup _{n\rightarrow\infty}\frac{d_n^{(3)}}{d_{n-1}^{(3)}}\leq \lim _{n\rightarrow\infty}[d_{n-1}^{(3)}]^{1/(n-1)}=\lambda _3(\widehat{F}). 
\end{eqnarray*}
Hence, by Lemmas \ref{Lemma2} and \ref{Lemma4}, we have that $\lambda _3(\widehat{F})=\lambda _2(\widehat{F})>\lambda _1(\widehat{F})$. Since $\lambda _1(\widehat{F^{-1}})=\lambda _3(\widehat{F})$ and $\lambda _2(\widehat{F^{-1}})=\lambda _2(\widehat{F})$, we have $\lambda _1(\widehat{F^{-1}})^2>\lambda _2(\widehat{F^{-1}})$. Applying Theorem 1 in \cite{dang-favre} and the fact that  the first degree sequence for $F^{-1}$ is the same as the third degree sequence for $F$, it follows easily that for all $j\geq 0$ we have
\begin{eqnarray*}
\lim _{n\rightarrow\infty}\frac{d_{n+j}^{(3)}}{d_n^{(3)}}=\lambda _3(\widehat{F})^j. 
\end{eqnarray*}

Assume that $d_n^{(3)}\leq 3/2 \cdot (d_{n-1}^{(3)}-d_{n-3}^{(3)})+2d_{n-2}^{(3)}$ for all $n\geq 5$. Dividing $d_{n-3}^{(3)}$ on both side of the inequality, and taking limit when $n\rightarrow\infty$, we obtain
\begin{eqnarray*}
\lambda _3(\widehat{F})^3\leq 3/2 \cdot (\lambda _3(\widehat{F})^2-1)+2\lambda _3(\widehat{F}). 
\end{eqnarray*}
Looking at the graph of the function $t\mapsto t^3-3/2 \cdot (t^2-1)+2t$ (for $t\in [1,\infty)$, recalling that dynamical degrees of a map are $\geq 1$), we see that the above inequality holds only if $\lambda _3(\widehat{F})\leq \zeta _2$. 
\end{proof}

Hence, one promising approach towards establishing $\lambda _2(\widehat{f}_c)\leq \zeta _2$ is to answer in the affirmative the following question: 

{\bf Question 1:} 

a) Is it true that $d_n^{(3)}\leq 3/2 \cdot (d_{n-1}^{(3)}-d_{n-3}^{(3)})+2d_{n-2}^{(3)}$ for $n$ large enough? 

b) Is it true that  $d_n^{(3)}\geq 2 d_{n-1}^{(3)}$ for $n$ large enough?

\subsubsection{Two approaches towards establishing the lower bound $\lambda _2(\widehat{f}_c)\geq \zeta _2$}

Again, to prove that $\lambda _2(\widehat{f}_c)\geq \zeta _2$ is the same as proving that $\lambda _3(\widehat{F})\geq \zeta _2$, and also is the same as proving that $\lambda _2(\widehat{F})\geq \zeta _2$. 

To this end, we have two approaches. One is again to base on the last two tables in the previous Subsection, while the other is based on a new viewpoint. This new viewpoint gives even more support to that we should have $\lambda _2(\widehat{f}_c)\geq \zeta _2$. 

\underline{\bf Approach 1:} As we mentioned, the last two tables seem to show that while $d_n^{(3)}$ (as well as $d_n^{(2)}$) does not satisfy the linear recurrence $w_n=3/2 \cdot (w_{n-1}-w_{n-3})+2w_{n-2}$, it is very close to satisfying the linear recurrence. More precisely, the difference is relatively small, in the sense that: 
\begin{eqnarray*}
\frac{d_n^{(3)}-3/2 \cdot (d_{n-1}^{(3)}-d_{n-3}^{(3)})-2d_{n-2}^{(3)}}{d_n^{(3)}}
\end{eqnarray*}
is small. This prompts us to ask the following question: 

{\bf Question 2:}

Is it true that 
\begin{eqnarray*}
\limsup _{n\rightarrow\infty}\frac{d_n^{(3)}-3/2 \cdot (d_{n-1}^{(3)}-d_{n-3}^{(3)})-2d_{n-2}^{(3)}}{d_n^{(3)}}\geq 0 ?
\end{eqnarray*}

We have the following result. 
\begin{lemma}
Assume that Question 2 has an affirmative answer. Then $\lambda _3(\widehat{F})\geq \zeta _2$. 

\label{Lemma6}\end{lemma}
\begin{proof}
Since $\lambda _1(\widehat{F})>1$, it follows by properties of dynamical degree {i.e. log-concavity} ) that also $\lambda _3(\widehat{F})>1$. Therefore, there cannot be a subsequence $n_k$ of the set of positive integers  so that $d_{n_k-1}^{(3)}-d_{n_k-3}^{(3)}\leq 0$ for all $k$. From this and the assumption that the answer to Question 2 is affirmative, we obtain 
\begin{eqnarray*}
\limsup _{n\rightarrow\infty}\frac{d_n^{(3)}}{d_{n-2}^{(3)}}\geq 2.  
\end{eqnarray*}
Then we have $\lambda _3(\widehat{F})\geq \sqrt{2}$. By Lemmas \ref{Lemma2} and \ref{Lemma4}, $2>\lambda _1(\widehat{F})$ and $\lambda _2(\widehat{F})=\max \{\lambda _1(\widehat{F},\lambda _3(\widehat{F})\}$. Therefore $\lambda _3(\widehat{F})^2>\lambda _2(\widehat{F})$.  As in the proof of Lemma \ref{Lemma5}, we then obtain $\lambda _3(\widehat{F})^3\geq 3/2 \cdot (\lambda _3(\widehat{F})^2-1)+2\lambda _3(\widehat{F})$. Since the polynomial $t^3-3/2 \cdot (t^2-1)-2t$ has 3 real roots with approximate values $-1.202$, $0.591$ and $\zeta _2\sim 2.1108$, while $\lambda _3(\widehat{F})\geq 1$ by definition, we conclude that we must have $\lambda _3(\widehat{F})\geq \zeta _2$. 

\end{proof}

Hence, a promising approach to establishing the lower bound $\lambda _3(\widehat{F})\geq \zeta _2$ is to solve in the affirmative Question 2. However, Approach 2 below seems to have more evidence to support than this Approach 1. 

\underline{\bf Approach 2:} In this approach, we compare the degree sequence $d_n^{(3)}$'s (respectively $d_n^{(2)}$'s) with the sequence $b_n$'s (correspondingly $c_n$'s). The sequence $b_n$'s satisfies the linear recurrence $w_n=3/2 \cdot (w_{n-1}-w_{n-3})+2w_{n-2}$ and has the first 3 initial values the same as that for the sequence $d_n^{(3)}$'s. Similarly, the sequence $c_n$'s satisfies the linear recurrence $w_n=3/2 \cdot (w_{n-1}-w_{n-3})+2w_{n-2}$ and has the first 3 initial values the same as that for the sequence $d_n^{(2)}$'s. From the experimental results in the previous Subsection, we get the following two tables. 

$$\begin{array}{c|c|c|c}
N & d_N^{(3)}=3^{\rm rd}~\hbox{degree of $F$} & b_N=3/2 \cdot (b_{N-1}-b_{N-3})+2b_{N-2} & \hbox{difference=$d_N^{(3)}-b_N$} \\ \hline
1 & 3 &3 & 0\\
2 & 7 &7 & 0\\
3 & 17 &17 &0 \\
4 & 37 & 35 & 2 \\ 
5 & 79 & 76 & 3\\
6 & 167 & 158.5&8.5 \\
7 & 353 &337.25  &15.75 \\
8 & 745 & 708.875 &  36.125\\ 
9 & 1571 &1500.0625  &70.9375 \\
10 & 3311 & 3161.96875 &149.03125\\
11 & 6977 & 6679.765625 & 297.234375\\
12 & 14701 &  14093.4921875 &607.5078125\\ 
13 & 30975 & 29756.81640625 & 1218.18359375\\
14&65263&62802.0560546875&2460.9439453125015\\
\hline
\end{array}$$
    
$$\begin{array}{c|c|c|c}
N & d_N^{(2)}=2^{\rm nd}~\hbox{degree of $F$} & c_N=3/2 \cdot (c_{N-1}-c_{N-3})+2c_{N-2} & \hbox{difference=$d_N^{(2)}-c_N$} \\ \hline
1 & 5 &5 & 0\\
2 & 9 &9 & 0\\
3 & 25 &25 &0 \\
3 & 49 & 48& 1\\
5 & 109 &108.5  &  0.5 \\ 
6 & 225 & 221.25 &3.75 \\
7 & 477 & 476.875&0.125 \\
8 & 1005 & 995.0625 & 9.9375\\
9 &  2117 &   2114.46875& 2.53125 \\ 
10 &  4465 &4446.515625  & 18.484375\\
11 &   9401 &9406.1171875&-5.1171875\\
12 &19817 & 19830.50390625 &-13.50390625 \\
13 & 41741 & 41888.216796875  &-147.216796875\\ 
14 & 87961&  88384.1572265625& -423.1572265625\\
\hline
\end{array}$$
 
 It seems very evident that we should have $d_n^{(3)}\geq b_n$ for all $n\geq 1$. This prompts us the following question. 
 
 {\bf Question 3:} Is it true that we have $d_n^{(3)}\geq b_n$ for all $n\geq 1$? (In the proof, we only need $d_n^{(3)}\geq \epsilon b_n$ for all $n\geq 1$ and a constant $\epsilon >0$.)
 
 We have the following result. 
 
 \begin{lemma} Assume that Question 3 has an affirmative answer. Then $\lambda _3(\widehat{f}_c)\geq \zeta _2$. 
\label{Lemma7}\end{lemma}
\begin{proof}
This follows easily from the fact that $\lim _{n\rightarrow\infty}b_n^{1/n}=\zeta _2$. 
\end{proof}

(Note that, on the other hand, it seems that for $n\geq 11$ then $d_n^{(2)}\leq c_n$. Using the latter inequality, we obtain only the upper bound $\lambda _2(F)\leq \zeta _2$, which is already discussed in the previous Subsubsection and not the lower bound wanted. However, one can ask whether $d_n^{(2)}\geq \epsilon c_n$ for all $n\geq 1$ and a constant $\epsilon >0$. This inequality seems to be supported by the data, and is also enough to deduce that $\lambda _3(\widehat{f}_c)\geq \zeta _2$.)

Hence, a promising approach towards showing $\lambda _2(\widehat{f}_c)\geq \zeta _2$ is to solve Question 3 in the affirmative. 

\subsubsection{An approach to part 3 of Conjecture \ref{Conjecture3}}

We divide this into two tasks: one concerning the Weaker estimate (which provides a counter-example to Conjecture \ref{Conjecture1}) and one concerning the Stronger estimate (which provides a counter-example to Conjecture \ref{Conjecture2}). 

\underline{\bf An approach towards the Weaker estimate in part 3 of Conjecture \ref{Conjecture3}:}

From the table for the number of isolated periodic points, we find that the sequence\\ $[\sharp IsoFix_{2n+1}(f_c)]^{(1/2n+1)}$ (for $0\leq n\leq 5$) is a {\bf decreasing sequence}: $4$,  $2.15443469003$, $2.13152551$, $2.10967780991$, $2.10744910267$, and $2.10703309279$. Also, we find that the fixed point set of $f_c^{2n+1}$ (for $0\leq n\leq 5$) consists of isolated points only. These facts do not look like a random coincidence, and hence naturally lead to the following question: 

{\bf Question 4:} Is it true that the sequence $[\sharp IsoFix_{2n+1}(f_c)]^{1/(2n+1)}$ (for $n=0,1,2,\ldots $) is a decreasing sequence? 

If Question 4 has an affirmative answer, then since $[\sharp IsoFix_{9}]^{1/9}=2.1074...<2.108$, we obtain right away a proof of the Weaker estimate in part 3 of Conjecture \ref{Conjecture3}. 

Here is a heuristic explanation for why Question 4 can have an affirmative answer: We know from Proposition \ref{TheoremMain} that the Lefschetz numbers $\{L(\widehat{f}_c^n)^{1/n}\}$ for the map $\widehat{f}_c:X_c\dashrightarrow X_c$ is decreasing. It seems that in this special case, we can localise this property to the map $f_c=\widehat{f}_c|_{Z_c}$. If this is so, and if we can show that the fixed point set of $f_c$ consists of only isolated points (as seen in the experiments), then Question 4 is solved in the affirmative (since then the Lefschetz number is the same as the number of fixed points).  

\underline{\bf An approach towards the Stronger estimate in part 3 of Conjecture \ref{Conjecture3}:}

This part is probably more difficult to establish. Our clue is that from the table we observe the following phenomenon: for all $0\leq n\leq 5$, we have $[\sharp IsoFix_{2n+1}(f_c)]^{1/{2n+1}}\geq [\sharp IsoFix_{2n+2}(f_c)]^{1/{2n+2}}$. Thus comes another question: 

{\bf Question 5:} Is it true that $[\sharp IsoFix_{2n+1}(f_c)]^{1/{2n+1}}\geq [\sharp IsoFix_{2n+2}(f_c)]^{1/{2n+2}}$ for all $n=0,1,2,\ldots $? 

If Question 5 has an affirmative answer, and if moreover the Weaker estimate in part 3 of Conjecture \ref{Conjecture3} holds, then the Stronger estimate in part 3 of Conjecture \ref{Conjecture3} follows readily.

Why Question 5 should have an affirmative answer could be again contributed to a localisation of Proposition \ref{TheoremMain}. While the log concavity of the sequence $IsoFix_n(f_c)$ is violated (see the next paragraph), still some parts could be preserved, allowing us to have affirmative answers to both Questions 4 and 5.  

In stead of Question 5, the following variant is also enough for our purpose 

{\bf Question 6:} Is it true that $\sharp IsoFix_{2n+1}(f_c)\geq \sharp IsoFix_{2n}(f_c)$ for all $n=1,2,3,\ldots $?

\subsubsection{A speculation on the exponential growth of the isolated periodic points}

Another curious phenomenon, which we do not need in resolving Conjecture \ref{Conjecture3}, is that  the sequence $[\sharp IsoFix_{2n+2}(f_c)]^{1/{2n+2}})$ (for $n=0,1,2,\ldots $) seems to be {\bf increasing}. We do not know of a possible explanation for this interesting phenomenon.

However, if this increasing phenomenon is true, and the phenomena mentioned in Questions 4 and 5 are also true, then we will obtain
\begin{eqnarray*}
\log 2.0890\leq \lim _{n\rightarrow\infty}\frac{\log \sharp IsoFix_{2n}(f_c)}{2n}\leq  \lim _{n\rightarrow\infty}\frac{\log \sharp IsoFix_{2n+1}(f_c)}{2n+1}\leq \log 2.1071,
\end{eqnarray*}
and it is then reasonable to speculate that indeed 
\begin{eqnarray*}
\lim _{n\rightarrow\infty}\frac{\log \sharp IsoFix_{n}(f_c)}{{n}}
\end{eqnarray*}
also exists. (In this case, the limit must be contained in the interval $[\log 2.0890,\log 2.1071]$.)

\section{Discussion} Our paper is the first one to display a potential counter-example to a major folklore conjecture, Conjecture \ref{Conjecture1}, in complex Dynamical Systems. At the same time, it presents explicitly the first time a potential example of a polynomial birational map on affine $4$-spaces whose dynamical degrees are algebraic numbers but not algebraic integers, thus is also a potential counter-example to another actively studied conjecture, that is Conjecture \ref{Conjecture2}. Even for the larger class of all dominant rational selfmaps of projective varieties, or compact K\"ahler manifolds, our maps are also the first potential counter-examples to the mentioned conjectures.   

The methods in previous work are not applicable to analyse our maps, at least for the moment. There are three common approaches toward computing dynamical degrees, which involve establishing a linear recurrence between the degrees of the iterates. The first one (for the first dynamical degree $\lambda _1$) is to observe directly a linear recurrence between the degrees of iterates of $f$. The second one is to construct a birational (or bimeromorphic) map $\pi :Y\rightarrow X$ for which the lifting map $f_Y$ satisfies a so-called algebraic stability under which $\lambda _j(f)$ is the spectral radius of the pullback map $f_Y^*:H^{j,j}(Y)\rightarrow H^{j,j}(Y)$. In the current literature, all of these approaches are only able to show that the concerned is a dynamical degree integer, which is not the case of our map. (A generalization of this is the recent work \cite{dang-favre} which shows that the pullback map on an associated Banach space $RZ(Y)$ is stable. However, unlike on the finite dimensional cohomology groups, there is no guarantee that the eigenvalues on $RZ(Y)$ are algebraic integers. See the next paragraph for more discussion on this method.)  Finally, the third one is to  use specialities of toric varieties and toric maps. A more recent method t in \cite{bell-etal}\cite{bell-etal2} is to use Diophantine approximation, but this can only prove that the concerned dynamical degree is a transcendental number, again not the case~for~our~maps. 

Now we discuss more why the results in \cite{dang-favre}, which is a very strong result on computing dynamical degrees, are not yet applicable to our maps. For example, the proof of Theorem~2 in \cite{dang-favre} - which asserts that for a proper polynomial map of the complex affine space $\mathbb{C}^d$ for which $\lambda _1^2>\lambda _2$, we always have that $\lambda _1$ is an algebraic number - cannot be used to establish part 2 of Conjecture \ref{Conjecture3} 
{for two reasons.}
First, while it seems true that $\lambda _3(\widehat{F})=\lambda _1(\widehat{F^{-1}})^2>\lambda _2(\widehat{F^{-1}})$ (for example, if one can show that $\lambda _3(\widehat{F})>\sqrt{2}$), the map $F^{-1}$ is not a polynomial map and hence the arguments used in Theorem 6.1 in \cite{dang-favre} - which rely on the valuation theory on the affine space - cannot be used. Second, even if the issue mentioned in the previous sentence can be resolved, it can be difficult to actually construct a matrix with rational coefficients whose spectral radius is exactly $\zeta _2$ and $\lambda _1(\widehat{F^{-1}})$ is also~an~eigenvalue.  

Previous work only concentrated on confirming Conjecture \ref{Conjecture1}, and hence the developed techniques cannot be used to prove Conjecture \ref{Conjecture3} which aims to disprove Conjecture \ref{Conjecture1}. 

Also, as mentioned, symbolic methods as Gr\"obner basis can be too slow and not having enough features needed to analyse our maps. 

To overcome the above disadvantages of the current techniques, we use the software Bertini instead of symbolic methods. Also, we observe that while it seems that the second degree sequence of $f_c$ does not satisfy any linear recurrence, there are some good comparisons (both lower and upper bounds) between the sequence and relevant linear recurrences. Since previous work has no mention on how to proceed if the number of periodic points of a map on a Zariski open set seems not as many as expected, we have to get around by looking whether the log-concavity phenomenon in Proposition \ref{TheoremMain} still holds for a specific map and a specific Zariski open set. 

\section{Conclusions} 
 
In this paper, we presented a simple family of birational maps on smooth affine quadric 3-folds, coming from a polynomial map on $\mathbb{C}^4$, which seem to be cohomologically hyperbolic while having less periodic points than  expected. Moreover, the second dynamical degree of these maps seem to be an algebraic number, but not an algebraic integer. 

This kind of maps requires the development of new tools/ideas stronger than those in the current literature, see the previous Discussion section. Among the theoretical tools, it needs the development of more effective non-proper intersections of varieties (not only for projetcive varieties, but also for affine varieties), which help to establish various log concavity phenomena.  It also requires taking into account periodic points located in proper analytic subvarieties. In dimension 2, only the case of bimeromorphic self-maps of surfaces has been dealt with in the literature \cite{saito}\cite{iwasaki-uehara}\cite{dinh-etal2}.  In higher dimensions, for holomorphic maps or more general meromorphic maps of compact K\"ahler manifolds satisfying certain technical conditions (including a so-called algebraic stability and an assumption on the sequence of currents of integration associated to the graphs of iterates of the map), there have been work in \cite{dinh-etal4} and \cite{vu}. Note that since the first dynamical degree of our map $\widehat{F^{-1}}$ is conjectured to be not an algebraic integer, there will be no birational model of its which is 1-algebraic~stable.

Besides Conjecture \ref{Conjecture2}, another consequence of Conjecture \ref{Conjecture1} is the following: 

\begin{conjecture} Let $X$ be a smooth complex projective variety, and  $f:X\dashrightarrow X$  a dominant rational map. If $f$ is cohomologically hyperbolic, then the set of isolated periodic points of $f$ is Zariski dense.   
\label{Conjecture5}\end{conjecture}
Conjecture \ref{Conjecture5} is independent of the birational model of $f$. Besides the cases where Conjecture \ref{Conjecture1} is known to be valid, one nontrivial case where Conjecture \ref{Conjecture5} is solved in the affirmative is when $X$ is of dimension $2$ and $f$ is a birational map \cite{xie}. (In this case, $f$ being cohomologically hyperbolic is equivalent to $\lambda _1(f)>1$, and the latter fact can be checked by a finite algorithm following the paper \cite{diller-favre} mentioned above. Indeed, there is a sufficient criterion for  $\lambda _1(f)>1$ in terms of some inequalities involving only several first terms in the degree sequence for $f$, see \cite{xie}.) In light of the map considered in Conjecture \ref{Conjecture3}, it is interesting to check whether Conjecture \ref{Conjecture5} holds for every cohomologically hyperbolic birational maps in dimension $3$.  
 
 Concerning the equality for the Gromov-Yomdin's and Dinh-Sibony's bounds on the topoligical entropy of a map $f$ in terms of its dynamical degrees, i.e., $h_{top}(f)\leq \max \log \lambda _j(f)$,  \cite{deThelin-vigny} shows that a generic birational map of $\mathbb{C}^d$ will obtain this bound. On the other hand, the proof therein uses the assumption that the maps concerned are algebraic stable (which is not satisfied for our map). The construction of the Banach space in \cite{dang-favre} makes use of the N\'eron-Severi groups of birational models of the given map, and for topological entropy a similar consideration using topological entropies of birational models of the given map was proposed in \cite{deThelin}. It is interesting to see whether the ideas in \cite{deThelin} work for our map $\widehat{f_c}$ (or~$\widehat{F}$). One can also try to check if the more general approach of using etale covers of the map (but then, must work with correspondences in general) in \cite{truong4} can work. 
 
 Among the experimental tools, one needs faster and less expensive (while still being highly reliable) methods to deal with numerical calculations for solving systems of equations, in particular those specially designed for compositions of simple maps. Indeed, the calculations presented in Section 2 are more or less at the upper limit of what current calculation methods can afford us (even for the very simple map $F$).   

Since dynamical degrees can be defined over any algebraically closed field \cite{truong}\cite{dang}, it is also meaningful to extend the study mentioned above to arbitrary algebraically closed fields.

\section{Acknowledgements and Declarations}

We thank Duc-Viet Vu for useful comments. C. B. and T. T. T. thank Research in pair grants and the hospitality from EPFL (Centre Bernoulli) and the University of Trento and FBK Foundation, as well as travel grants from Trond Mohn Foundation, for their supports. They are also supported partly by Young Research talents grant 300814 from Research Council of Norway. {The first author C.B. was also partially supported by GNSAGA of INdAM and by PRIN Variet\'a reali e complesse: geometria, topologia e
analisi armonica.}
J. D. H. was supported partly by National Science Foundation grants CCF 1812746
and CCF 2331400.  

{\bf Declarations of interests:} none.

\end{document}